\input ./style/arxiv-vmsta.cfg
\documentclass[numbers,compress,v1.0.1]{vmsta}
\usepackage{dcolumn}
\usepackage{amsmath}

\volume{2}
\issue{2}
\pubyear{2015}
\firstpage{165}
\lastpage{171}
\doi{10.15559/15-VMSTA29}% Updated by VTEXPTS2LaTeX.exe, 27.07.2015
\startlocaldefs
\newtheorem{thm}{Theorem}

\theoremstyle{definition}
\newtheorem{remark}{Remark}

\newcolumntype{d}[1]{D{.}{.}{#1}}

\allowdisplaybreaks
\endlocaldefs

\begin{document}
\begin{frontmatter}

\title{Fast $L_2$-approximation of integral-type functionals of~Markov
processes}
\author{\inits{Iu.}\fnm{Iurii}\snm{Ganychenko}}\email
{iurii\_ganychenko@ukr.net}%\corref{cor1}
%\cortext[cor1]{Corresponding author.}

\address{Taras Shevchenko National University of Kyiv, Kyiv, Ukraine}

\markboth{Iu. Ganychenko}{Fast $L_2$-approximation of integral-type
functionals of Markov processes}

\begin{abstract}
In this paper, we provide strong $L_2$-rates of approximation of the
integral-type functionals of Markov processes by integral sums. We
improve the method developed in \cite{kul-gan}. Under assumptions on
the process formulated only in terms of its transition probability
density, we get the accuracy that coincides with that obtained in \cite
{Kohatsu-Higa} for a one-dimensional diffusion process.
\end{abstract}

\received{8 July 2015}% Updated by VTEXPTS2LaTeX.exe, 27.07.2015 13:09
\revised{20 July 2015}% Updated by VTEXPTS2LaTeX.exe, 27.07.2015 13:09
\accepted{22 July 2015}% Updated by VTEXPTS2LaTeX.exe, 27.07.2015 13:09

\publishedonline{28 July 2015}
\end{frontmatter}
\section{Introduction}

Let $X_t,\, t\geq0$, be a Markov process with values in $\mathbb
{R}^d$. Consider the following objects:
\begin{enumerate}
\item[1)] the integral functional
\[
I_T(h)=\int_0^Th(X_t)
\, dt
\]
of this process;
\item[2)] the sequence of integral sums
\[
I_{T,n}(h)={T\over n}\sum_{k=0}^{n-1}h(X_{(kT)/n}),
\quad n\geq1.
\]
\end{enumerate}

In this paper, we establish \emph{strong $L_2$-approximation rates},
that is, the bounds for
\[
E \bigl|I_{T}(h)-I_{T,n}(h) \bigr|^2.
\]
The current research is mainly motivated by the recent papers \cite
{kul-gan} and \cite{Kohatsu-Higa}.

In \cite{Kohatsu-Higa}, strong $L_p$-approximation rates are considered
for an important particular case where $X$ is a one-dimensional
diffusion. The approach developed in this paper contains both the
Malliavin calculus tools and the Gaussian bounds for the transition
probability density of the process $X$, and relies substantially on the
structure of the process.

Another approach to that problem has been developed in \cite{kul-gan}.
This approach is, in a sense, a modification of Dynkin's theory of
continuous additive functionals (see~\cite{Dynkin}, Chap.~6) and also
involves the technique similar to that used in the proof of the
classical Khasminskii lemma (see, e.g., \cite[Lemma~2.1]{Sznitman}).
This approach allows us to obtain strong $L_p$-approximation rates
under assumptions on the process $X$ formulated only in terms of its
transition probability density.

For a bounded function $h$, the strong $L_p$-rates of approximation of
the integral functional $I_T(h)$ obtained in \cite{kul-gan} essentially
coincide with those established in \cite{Kohatsu-Higa}. However, under
additional regularity assumptions on the function~$h$ (e.g., when $h$
is H\"{o}lder continuous), the rates obtained in \cite{Kohatsu-Higa}
are sharper (see \cite[Thm.~2.2]{kul-gan} and \cite[Thm.~2.3]{Kohatsu-Higa}).

In this note, we improve the method developed in \cite{kul-gan}, so
that under the assumption of the H\"{o}lder continuity of $h$, the
strong $L_2$-approximation rates coincide with those obtained in \cite
{Kohatsu-Higa}, preserving at the same time the advantage of the method
that the assumptions on the process $X$ are quite general and do not
essentially rely on the structure of the process.

\section{Main result}

In what follows, $P_x$ denotes the law of the Markov process $X$
conditioned by $X_0=x$, and $\mathbb{E}_x$ denotes the expectation with
respect to this law. Both the absolute value of a real number and the
Euclidean norm in $\mathbb{R}^d$ are denoted by~$|\cdot|$.

We make the following assumption on the process $X$.

\medskip
\noindent\textbf{A.} The process $X$ possesses a transition probability density
$p_t(x,y)$ that is differentiable with respect to $t$ and satisfies the
following estimates:
\begin{equation}
\label{dens_bound} p_t(x,y)\leq C_Tt^{-d/\alpha} Q
\bigl(t^{-1/\alpha}(x-y) \bigr), \quad t\leq T,
\end{equation}
\begin{equation}
\label{der_bound}
\bigl|\partial_tp_t(x,y) \bigr|
\leq
C_Tt^{-1-d/\alpha} Q \bigl(t^{-1/\alpha}(x-y) \bigr), \quad t\leq
T,
\end{equation}
\begin{equation}
\label{der2_bound}
\bigl|\partial^2_{tt}p_t(x,y) \bigr|
\leq
C_Tt^{-2-d/\alpha} Q \bigl(t^{-1/\alpha}(x-y) \bigr), \quad t\leq T,
\end{equation}
for some fixed \mbox{$\alpha\in(0,2]$} and some distribution density $Q$ such
that\break $\int_{\mathbb{R}^d} |z|^{2\gamma} Q(z)\, dz < \infty$. Without
loss of generality, we assume that in \eqref{dens_bound}--\eqref
{der2_bound} $C_T\geq1$.

We assume that the function $h$ satisfies the H\"{o}lder condition with
exponent $\gamma\in(0, \alpha/2]$, that is,
\[
\|h\|_\gamma:=\sup_{x\not=y}{|h(x)-h(y)|\over|x-y|^\gamma}<
\infty.
\]

Now we formulate the main result of the paper.

\begin{thm}\label{mr}
Suppose that Assumption \textbf{A} holds. Then
\begin{equation*}
\mathbb{E}_x
\bigl|I_{T}(h)-I_{T,n}(h) \bigr|^2
\leq%
\begin{cases}
D_{T,\gamma,\alpha,Q} C_{\gamma,\alpha}\|h\|^2_\gamma
n^{-(1+2\gamma
/\alpha)},& \gamma\neq\alpha/2, \\
D_{T,\gamma,\alpha,Q} \|h\|^2_\gamma n^{-2} \ln n, & \gamma= \alpha/2,
\end{cases} %
\
\end{equation*}
where
\begin{align*}
&D_{T,\gamma,\alpha,Q} = 8C^2_T
T^{2+2\gamma/\alpha} \int_{\mathbb
{R}^d} |z|^{2\gamma} Q(z)\, dz,
\\
& C_{\gamma,\alpha} = \max \biggl\{ (1-2\gamma/\alpha )^{-1}(2\gamma/
\alpha)^{-1}, \, \max_{n\geq1} \biggl( \frac{(\ln
n)^2}{n^{1-2\gamma/\alpha}}
\biggr) \biggr\}.
\end{align*} %
\end{thm}

We provide the proof of Theorem \ref{mr} in Section \ref{proof}.

\begin{remark}
Any diffusion process satisfies conditions {(\ref
{dens_bound})--(\ref{der2_bound})} with $\alpha= 2$,
$Q(x)=c_1e^{-c_2|x|^2}$, and properly chosen $c_1, c_2$ {(}see
{\cite{kul-gan})}. In the case where $X$ is a one-dimen\-sional
diffusion, Theorem \ref{mr} provides the same rates of convergence as
those obtained in {\cite{Kohatsu-Higa}} {(}see Theorem 2.3 in
{\cite{Kohatsu-Higa})}.
\end{remark}

\begin{remark}
Similarly to {\cite{kul-gan}}, we formulate the assumption on the
process $X$ only in terms of its transition probability density.
Condition \textbf{A}, compared with condition \textbf{X} {(}cf.
{\cite{kul-gan})}, contains the additional assumption {(\ref
{der2_bound})}.
\end{remark}

\section{Proof of Theorem \ref{mr}}
\label{proof}
\begin{proof}\begingroup%\abovedisplayskip=8pt\belowdisplayskip=8pt
For $t\in[kT/n, (k+1)T/n)$, denote
\[
\eta_n(t)={kT\over n}, \qquad\zeta_n(t)=
{(k+1)T\over n},
\]
and put $\Delta_n(s):=h(X_s) - h(X_{\eta_n (s)})$, $s \in[0,T]$.

By the Markov property of $X$, for any $r<s$, we have
\begin{align*}
\mathbb{E}_x|X_s-X_r|^{2\gamma}
&{}=\mathbb{E}_x \int_{\mathbb{R}^d}p_{s-r}(X_r,z)|X_r-z|^{2\gamma}\,dz
\\
&{}\leq C_T \mathbb{E}_x \int_{\mathbb{R}^d}(s-r)^{-d/\alpha}Q
\bigl((s-r)^{-1/\alpha}(X_r-z) \bigr)|X_r-z|^{2\gamma}
\,dz
\\
&{} = C_T(s-r)^{2\gamma/\alpha}\int_{\mathbb{R}^d}|z|^{2\gamma}
Q(z)\,dz.
\end{align*} %

Therefore, using the inequality $s-\eta_n(s) \leq T/n$, $s \in[0,T]$ and
the H\"{o}lder continuity of the function $h$, we obtain:
\begin{equation}
\label{Lip_bound} \mathbb{E}_x\bigl|\Delta_n(s)\bigr|^2
\leq C_T T^{2\gamma/\alpha} \biggl( \int_{\mathbb{R}^d}
|z|^{2\gamma} Q(z)\, dz \biggr) \|h\|^2_{\gamma
}n^{-2\gamma/\alpha}.
\end{equation}

Split
\begin{equation}
\label{initial} \mathbb{E}_x\bigl|I_{T}(h)-I_{T,n}(h)\bigr|^2=2
\mathbb{E}_x\int_{0}^{T}\int
_{s}^{T}\Delta_n(s)
\Delta_n(t)\, dt\,ds = J_1+J_2+J_3,
\end{equation}
where
\begin{align*}
&J_1= 2\mathbb{E}_x \int_{0}^{T} \int_{s}^{\zeta_n(s)+T/n}
\Delta_n(s) \Delta_n(t)\,dt\,ds,
\\[3pt]
& J_2= 2\mathbb{E}_x \int_{0}^{T/n}
\int_{\zeta_n(s)+T/n}^{T}\Delta_n(s)
\Delta_n(t)\,dt\,ds,
\\[3pt]
& J_3 = 2\mathbb{E}_x \int_{T/n}^{T}
\int_{\zeta_n(s)+T/n}^{T}\Delta_n(s)
\Delta_n(t)\, dt\,ds.
\end{align*} %
For $|J_1|$ and $|J_2|$, the estimates can be obtained in the same way.
Indeed, using the Cauchy inequality and (\ref{Lip_bound}), we get
\begin{align*}
|J_1|
&{}\leq2\int_{0}^{T}
\int_{s}^{\zeta_n(s)+T/n} \bigl(\mathbb{E}_x\bigl|\Delta_n(s)\bigr|^2 \bigr)^{1/2}
\bigl(\mathbb{E}_x\bigl|\Delta_n(t)\bigr|^2\bigr)^{1/2}\,dt\,ds
\\[3pt]
&{}\leq2C_T T^{2\gamma/\alpha}\|h\|^2_\gamma
\biggl( \int_{\mathbb{R}^d} |z|^{2\gamma} Q(z)\, dz \biggr)
n^{-2\gamma/\alpha} \int_{0}^{T} \bigl(T/n+\zeta_n(s)-s \bigr) \, ds
\\[3pt]
&{}\leq4C_T T^{2+2\gamma/\alpha}\|h\|^2_\gamma
\biggl( \int_{\mathbb{R}^d} |z|^{2\gamma} Q(z)\, dz \biggr)
n^{-(1+ 2\gamma/\alpha)}.
\end{align*} %
In the last inequality, we have used the inequality $\zeta_n(s)-s \leq
T/n$, $s \in[0,T]$.
Similarly,
\[
|J_2| \leq2C_T T^{2+2\gamma/\alpha}\|h\|^2_\gamma
\biggl(\int_{\mathbb
{R}^d} |z|^{2\gamma} Q(z)\, dz \biggr)
n^{-(1+ 2\gamma/\alpha)}.
\]
Now we proceed to the estimation of $|J_3|$, which is the main part of
the proof.
Observe that the following identities hold:
\begin{align}
\int_{\mathbb{R}^d}\partial^2_{uv}p_u(x,y)p_{v-u}(y,z)\,dz
&{} = \partial^2_{uv} p_u(x,y) \int
_{\mathbb{R}^d} p_{v-u}(y,z)\, dz\nonumber\\[3pt]
&{} = \partial^2_{uv} p_u(x,y) = 0, \quad y \in \mathbb{R}^d, \label{id1}\\[3pt]
%\end{align}%
%
%\begin{align}
\int_{\mathbb{R}^d} \partial^2_{uv} p_u(x,y) p_{v-u}(y,z) \, dy
&{} = \partial^2_{uv}\int_{\mathbb{R}^d}
p_u(x,y) p_{v-u}(y,z)\, dy
\nonumber\\[3pt]
&{} = \partial^2_{uv} p_v(x,z) = 0, \quad z \in
\mathbb{R}^d, \label{id2}
\end{align} %
where in (\ref{id1}) we used that $\int_{\mathbb{R}^d} p_{r}(y,z)\, dz
= 1$,  $r>0$, $y \in\mathbb{R}^d$, and in (\ref{id2}) we used the
Chapman--Kolmogorov equation.

We have:
\begin{align}
J_3
&{} = 2 \int_{T/n}^{T}\int_{\zeta_n(s)+T/n}^{T}\int_{\mathbb{R}^d} \int_{\mathbb{R}^d} h(y)h(z)
\bigl[p_s(x,y)p_{t-s}(y,z)
\nonumber\\
&\quad{} - p_{\eta_n(s)}(x,y)p_{t-\eta_n(s)}(y,z) - p_{s}(x,y)p_{\eta
_n(t)-s}(y,z)
\nonumber\\
&\quad{} + p_{\eta_n(s)}(x,y)p_{\eta_n(t)-\eta_n(s)}(y,z)\bigr]\,dz\,dy\,dt\,ds
\nonumber\\
&{} = 2 \int_{T/n}^{T}\int_{\zeta_n(s)+T/n}^{T}\int_{\mathbb{R}^d} \int_{\mathbb{R}^d} \int_{\eta_n(s)}^{s} \int_{\eta_n(t)}^{t}
h(y)h(z) \partial^2_{uv} \bigl(p_u(x,y)\nonumber\\
&\quad{} \times p_{v-u}(y,z)\bigr)\,dv\,du\,dz\,dy\,dt\,ds
\nonumber\\[2pt]
&{}= - \int_{T/n}^{T}\int_{\zeta_n(s)+T/n}^{T}
\int_{\mathbb{R}^d} \int_{\mathbb{R}^d} \int_{\eta_n(s)}^{s} \int_{\eta_n(t)}^{t}
\bigl(h(y)-h(z)\bigr)^2 \partial^2_{uv}\bigl(p_u(x,y)
\nonumber\\
&\quad{} \times p_{v-u}(y,z)\bigr)\, dv\,du\,dz\,dy\,dt\,ds,\label{j3}
\end{align} %
where in the last identity we have used (\ref{id1}) and
(\ref{id2}).\endgroup

Further, we have
\[
\partial^2_{uv} p_u(x,y) p_{v-u}(y,z)=
p_u(x,y)\partial^2_{rr}p_r(y,z)
\big|_{r=v-u} + \partial_{u}p_u(x,y)
\partial_{r}p_r(y,z) \big|_{r=v-u}.
\]
Then, using condition \textbf{A} and the H\"{o}lder continuity of the
function $h$, we obtain
\begin{align}
&\int_{\mathbb{R}^d} \int
_{\mathbb{R}^d} \bigl(h(y)-h(z)\bigr)^2\vert \partial
^2_{uv} \bigl(p_u(x,y) p_{v-u}(y,z)
\bigr) \vert\, dz\,dy\nonumber
\\
&\quad{} \leq C_T^2 \|h\|^2_{\gamma} \biggl(
\int_{\mathbb{R}^d} |z|^{2\gamma} Q(z)\, dz \biggr) \bigl(
(v-u)^{2\gamma/\alpha-2} + (v-u)^{2\gamma
/\alpha-1} u^{-1} \bigr).
\label{hh}
\end{align} %

Therefore, according to (\ref{j3}) and (\ref{hh}),
\begin{align}
&{}|J_3| \leq
C_T^2 \|h\|^2_{\gamma} \biggl( \int
_{\mathbb{R}^d} |z|^{2\gamma} Q(z)\, dz \biggr)\nonumber
\\
&\quad{}\times\int_{T/n}^{T}\int_{\zeta_n(s)+T/n}^{T}
\int_{\eta_n(s)}^{s} \int_{\eta_n(t)}^{t}
\bigl( (v-u)^{2\gamma/\alpha-2} + (v-u)^{2\gamma
/\alpha
-1} u^{-1} \bigr) \,
dv\,du\,dt\,ds.\label{j31}
\end{align}

Denote $a_{\alpha, \gamma} (u,v): = (v-u)^{2\gamma/\alpha-2} +
(v-u)^{2\gamma/\alpha-1} u^{-1}$. Then
\begin{align*}
&\int_{T/n}^{T}\int_{\zeta_n(s)+T/n}^{T} \int_{\eta_n(s)}^{s}
\int_{\eta_n(t)}^{t} a_{\alpha, \gamma} (u,v) \, dv\,du\,dt\,ds
\\[2pt]
&\quad{} =
\sum\limits_{i=1}^{n-1}
\sum\limits_{j=i+2}^{n-1}
\int_{iT/n}^{(i+1)T/n}
\int_{jT/n}^{(j+1)T/n}
\int_{iT/n}^{s}\int_{jT/n}^{t}a_{\alpha, \gamma} (u,v) \,
dv\,du\,dt\,ds
\\[2pt]
&\quad{} =
\sum\limits_{i=1}^{n-1}
\sum\limits_{j=i+2}^{n-1}
\int_{iT/n}^{(i+1)T/n}
\int_{jT/n}^{(j+1)T/n}
\int_{u}^{(i+1)T/n}
\int_{v}^{(j+1)T/n}a_{\alpha, \gamma} (u,v) \,
dt\,ds\,dv\,du
\\[2pt]
&\quad{}\leq T^2 n^{-2}
\sum\limits_{i=1}^{n-1}
\sum\limits_{j=i+2}^{n-1}
\int_{iT/n}^{(i+1)T/n}
\int_{jT/n}^{(j+1)T/n}
a_{\alpha, \gamma} (u,v) \,
dv\,du
\\[2pt]
&\quad{} = T^2 n^{-2}
\sum\limits_{i=1}^{n-1}
\int_{iT/n}^{(i+1)T/n}
\int_{(i+2)T/n}^{T}a_{\alpha, \gamma} (u,v) \,
dv\,du,
\end{align*} %
where in the fourth line we used that, for $u \in[iT/n,(i+1)T/n)$ and
$ v \in[jT/n,\allowbreak (j+1)T/n)$, we always have
$(i+1)T/n - u \leq T/n$ and $(j+1)T/n - v \leq T/n$.

Thus, from (\ref{j31}) we obtain
\begin{equation}
\label{j32ss} |J_3| \leq C_T^2
T^2 \|h\|^2_{\gamma} \biggl( \int
_{\mathbb{R}^d} |z|^{2\gamma} Q(z)\, dz \biggr) n^{-2}
(S_1+S_2),
\end{equation}
where
\begin{align*}
&S_1 = \sum\limits_{i=1}^{n-1} \int_{iT/n}^{(i+1)T/n}
\int_{(i+1)T/n}^{T} (v-u)^{2\gamma/\alpha-2} \, dv\,du,
\\
& S_2 = \sum\limits_{i=1}^{n-1}
\int_{iT/n}^{(i+1)T/n}
\int_{(i+2)T/n}^{T}(v-u)^{2\gamma/\alpha-1} u^{-1} \, dv\,du.
\end{align*} %

We estimate each term separately. In what follows, we consider the case
$\gamma< \alpha/2$; the case of $\gamma= \alpha/2$ is similar and
therefore omitted.
We have
\begin{align}
S_1
&{}\leq(1-2\gamma/\alpha)^{-1}
\sum\limits_{i=1}^{n-1}
\int_{iT/n}^{(i+1)T/n}
\bigl((i+1)T/n-u\bigr)^{2\gamma/\alpha-1} \,du\nonumber
\\
&{}=(1-2\gamma/\alpha)^{-1} (2\gamma/\alpha)^{-1}
\sum\limits_{i=1}^{n-1}
\bigl((i+1)T/n-iT/n\bigr)^{2\gamma/\alpha}
\nonumber
\\
&{}\leq(1-2\gamma/\alpha)^{-1} (2\gamma/\alpha)^{-1}
T^{2\gamma/\alpha} n^{1-2\gamma/\alpha}
\leq
C_{\gamma,\alpha} T^{2\gamma/\alpha}n^{1-2\gamma/\alpha}.\label{s1}
\end{align} %

Finally, since $v-u \leq T$ for $0\leq u < v \leq T$, we have
\begin{align}
S_2
&{}\leq T^{2\gamma/\alpha}
\sum\limits_{i=1}^{n-1}
\int_{iT/n}^{(i+1)T/n}
\int_{(i+2)T/n}^{T} (v-u)^{-1} u^{-1}
\,dv\,du
\nonumber\\
&{}\leq T^{2\gamma/\alpha}
\sum\limits_{i=1}^{n-1}
\Biggl(\int_{iT/n}^{(i+1)T/n} u^{-1} \, du \Biggr)
\Biggl(\int_{(i+2)T/n}^{T} \bigl(v-(i+1)T/n\bigr)^{-1}\,dv\Biggr)
\nonumber\\
&{} \leq T^{2\gamma/\alpha} \ln n
\sum\limits_{i=1}^{n-1}
\Biggl(\int_{iT/n}^{(i+1)T/n} u^{-1} \, du \Biggr) = T^{2\gamma/\alpha} (\ln n)^2
\nonumber\\
&{} \leq C_{\gamma,\alpha} T^{2\gamma/\alpha} n^{1-2\gamma/\alpha}.
\label{s2}
\end{align}

Combining inequality (\ref{j32ss}) with (\ref{s1}) and (\ref{s2}),
we derive
\[
|J_3| \leq2C_{\gamma,\alpha} C_T^2
T^{2+2\gamma/\alpha} \|h\|^2_{\gamma
} \biggl( \int
_{\mathbb{R}^d} |z|^{2\gamma} Q(z)\, dz \biggr) n^{-(1+2\gamma
/\alpha)}.\qedhere
\]
\end{proof}

\end{document}